\documentclass[journal]{IEEEtran}

\usepackage[cmex10]{amsmath}
\usepackage{amssymb}
\usepackage{xy}
\usepackage{url}
\usepackage[dvips]{graphicx}
\DeclareGraphicsExtensions{.eps}

\newcommand{\f}{\mathbf}
\newcommand{\g}{\mathfrak}

\newfont{\hiera}{cmsy10 scaled 2488} 
\newfont{\hierb}{cmsy10 scaled 1728}
\newfont{\hierc}{cmsy10 scaled 1200}

\newcommand{\Bigast}{
\mathop{\vphantom{\sum}\lower2.5pt\hbox{\hiera\char3}}}%

\newcommand{\Bigtimes}{
\mathop{\vphantom{\sum}\lower2.5pt\hbox{\hiera\char2}}}%

\newtheorem{theorem}{Theorem}

\begin{document}
\title{An Extension to an Algebraic Method for Linear Time-Invariant System and Network Theory:\\ {The full AC-Calculus}\\
{\huge (2007 Re-Release)}}

\author{Eberhard H.-A.\ Gerbracht
\thanks{This article first appeared in: The SMACD Committee (Ed.): Proceedings of the Fifth International Workshop on Symbolic Methods and Applications in Circuit Design, SMACD '98, Kaisers\-lautern, October 8-9, 1998. Kaisers\-lautern, 1998, pp.~134--139. Due to the low distribution of these proceedings, the author has decided to make the article available to a wider audience through the arXiv.}
\thanks{At the time of origin of this paper the author was with the Institut f\"ur Netz\-werk\-theorie und Schaltungstechnik, Technische Universit\"at Braunschweig, D-38106 Braunschweig, Germany.}
\thanks{His current (09/17/07) address is Bismarckstra\ss e 20, D-38518 Gifhorn, Germany. Current e-mail: \tt{e.gerbracht@web.de}}}
 
\maketitle

\begin{abstract}
Being inspired by phasor analysis in linear circuit theory, and its algebraic counterpart -- the AC-(operational)-calculus for sinusoids developed by W.~Marten and W.~Mathis -- we define a complex structure on several spaces of real-valued elementary functions. This is used to algebraize inhomogeneous linear ordinary differential equations with inhomogenities stemming from these spaces. Thus we deduce an effective method to calculate particular solutions of these ODEs in a purely algebraic way. 
\end{abstract}

{\small
\smallskip

}

\section{Introduction}
When considering linear time-invariant systems, i.e.\ (at least in this paper) when considering linear ordinary differential equations with constant coefficients, a number of ad-hoc-methods are used to calculate the response to certain kind of input functions. In this regard, phasor analysis, the AC-calculus for sinusoids developed by W.\ Marten and W.\ Mathis, and \lq\lq Ans\"atze\rq\rq\ according to the special form of the input come to mind. 

We will show that all of these methods fit under a common heading -- a heading, which we call the \lq\lq full AC-calculus\rq\rq. This calculus can be characterized on one hand by its tendency to group certain functions into linear spaces of functions. Thus we are able to use the full force of linear algebra. Another characteristic of this calculus is the extensive use of complex structures, and the quite astonishing fact that a number of spaces of real-valued functions allow the introduction of a multiplication with complex scalars, so that they carry the additional structure of a complex linear space. This quality allows for nice formulae and elegant algorithms.

This paper is still tentative, with proofs having been left out. Furthermore we do not yet show all possibilities of the full AC-calculus. Nevertheless the results are sound, and we hope to give our readers a glimpse at what might be possible. 

\section{Notation and Review of known facts}
Let $I$ be the interval $[0,\infty).$ The set of $n$-times continuously differentiable real-valued functions on the interval $I$ will be denoted by $C^{n}(I)$. We remind the reader of the fact that $C^{n}(I)$ is an infinite dimensional real linear space, i.e.\ superposition holds in $C^{n}(I)$.

Given $x\in C^n(I)$ we define a linear operator $L$ on $C^n(I)$ by setting
\begin{equation}
\label{Lop}
L(x) := a_n x^{(n)} + a_{n-1} x^{(n-1)} + \dots + a_1 \dot x + a_0 x
\end{equation}
with $a_0,\dots,a_n\in \mathbf R$. We call the operator $L$ {\em normalized} if $a_n = 1$. To $L$ we assign the {\em characteristic polynomial} $p_L$, which is
\begin{equation}
\label{CharPol}
p_L(s) = a_n s^n + a_{n-1}s^{n-1}+ \dots+a_1 s + a_0.
\end{equation}
A {\em inhomogeneous linear (ordinary) differential equation} can now be written very conveniently as
\begin{equation}
\label{inhomogeneDgl}
L(x) = r,
\end{equation}
where $r$ is a continous function on $I$. The corresponding {\em homogeneous linear ODE} is given by
\begin{equation}
\label{homogeneDgl}
L(x) = 0.
\end{equation}
Furthermore for $x\in C^n(I)$ and $t_0\in I$ we define the vector ${\mathbf x}(t_0)$ as
\begin{equation}
\label{Anfangsbed}
{\mathbf x}(t_0) := (x(t_0), \dot x(t_0),\dots, x^{(n-1)}(t_0))^\top.
\end{equation}
A typical {\em initial-value problem} can be written in the form
\begin{equation}
\label{IVP}
L(x)=r,\quad {\mathbf x}(0) = {\mathbf x}_0
\end{equation}
where $r$ is like above and ${\mathbf x}_0\in{\mathbf R}^n$.

The following theorem (cp.\ \cite{WalterDgl}, \S 20, Theorem 1, \S 19 I. and V.) summarizes the known facts about the solutions of (\ref{homogeneDgl}), (\ref{inhomogeneDgl}) and (\ref{IVP}).

\bigskip
\theorem
\label{Walter} 
The solutions of (\ref{homogeneDgl}) form an $n-$dimensional $\mathbf R$-linear subspace $X$ of $C^n(I)$. If $\lambda$ is a real zero of the characteristic polynomial with multiplicity $k$, then the functions
\begin{equation}
\label{Base1}
e^{\lambda t}, te^{\lambda t},\dots,t^{k-1}e^{\lambda t}
\end{equation}
form a $k-dimensional$ subspace of $X$. If $\lambda + j\omega, \lambda,\omega\in {\mathbf R}$, is a complex zero of the characteristic polynomial with multiplicity $k$ -- thus, because the coefficients $a_0,\dots,a_n$ are real, $\lambda - j\omega$ is a complex zero of multiplicity $k$, too -- then the functions
\begin{align}
\label{Base2a}
&e^{\lambda t}\sin(\omega t), te^{\lambda t}\sin(\omega t),\dots,t^{k-1}e^{\lambda t}\sin(\omega t),\\
\label{Base2b}
&e^{\lambda t}\cos(\omega t), te^{\lambda t}\cos(\omega t),\dots,t^{k-1}e^{\lambda t}\cos(\omega t)
\end{align} 
form a $2k-dimensional$ subspace of $X$. Consequently, the zeroes of the characteristic polynomial $p_L$ counted with multiplicity determine $n$ $\mathbf R$-linear independent solutions of (\ref{homogeneDgl}).

The solutions of the inhomogeneous linear ODE (\ref{inhomogeneDgl}) are given by
\begin{equation}
x = x_p + x_h
\end{equation}
where $x_p$ is a fixed solution of (\ref{inhomogeneDgl}) and $x_h$ is any solution of the homogeneous linear ODE (\ref{homogeneDgl}). Thus the solutions of (\ref{inhomogeneDgl}) form an {\em affine space} $A$, which is given by
\begin{equation}
A = x_p + X.
\end{equation}
Finally the initial-value problem (\ref{IVP}) uniquely determines a "point" $x_0\in A$, i.e.\ a function $x_0\in C^n(I)$, which then can be written in the form
\begin{equation}
\label{Sum}
x_0 = x_p + x_{h0}.
\end{equation}
with a suitable function $x_{h0}\in X$.
\endproof

\smallskip
In the language of  system theory we refer to $r$ as the {\em input}, to the operator $L$ as a {\em linear time-invariant system} and to the solution $x_0$ of the initial-value problem (\ref{IVP}) as the {\em output} or {\em response}.

\smallskip
Now we can formulate the time-honoured general strategy to solve the initial-value problem (\ref{IVP}):
\begin{itemize}
\item[\rm 1)] \label{step1} Find a particular solution $x_p$ to the inhomogeneous linear ODE (\ref{inhomogeneDgl}).
\item[\rm 2)] \label{step2} Determine the zeroes of the characteristic polynomial $p_L$ with multiplicities, and calculate the $\mathbf R$-linear independant solutions $x_{h1},\dots,x_{hn}$ to (\ref{homogeneDgl}) according to (\ref{Base1}), (\ref{Base2a}) and (\ref{Base2b}).
\item[\rm 3)] \label{step3} Make an "Ansatz"
\begin{equation}
x := x_p + \sum_{i=1}^n b_i x_{hi},
\end{equation}
determine the derivatives $\dot x, \dots, x^{(n-1)}, x^{(n)}$ and solve the system of inhomogeneous linear equations resulting from setting ${\mathbf x}(0) = {\mathbf x}_0$ with respect to the $b_i$. 
\end{itemize}

\bigskip
Steps 2) and 3), while sometimes necessitating annoying and not very easy to do calculations, conceptually do not pose any problems. The main difficulty of the above algorithm is the task set up in 1): finding at least one particular solution to a given inhomogeneous linear ODE.

Indeed in theoretical electrical engineering, there are several very well known approaches to this problem. All of them stem from the fact that, although the solution to the initial-value problem (\ref{IVP}) is a unique function, its decomposition into a sum of two functions as in (\ref{Sum}) is by far not unique (cp.\ \cite{DesoerKuh}, Chap.\ 4.3, in particular the introductory paragraph of section 4.3.2). We will substantiate this remark by two examples.

\subsection{Zero-State and Zero-Input Response (\cite{DesoerKuh}, Chap.\ 6.1)}
In the situation given by (\ref{IVP}) the unique solution of
\begin{equation}
\label{zeroinput}
L(x_{h0})=0,\quad {\mathbf x}_{h0}(0) = {\mathbf x}_0
\end{equation}
is called the {\em zero-input solution/response} to the initial-value problem,
the unique solution of
\begin{equation}
\label{zerostate}
L(x_{p0})=r,\quad {\mathbf x}_{p0}(0) = {\mathbf 0}
\end{equation}
is called the {\em zero-state (solution/response)}. Clearly, $x_{h0}$ is a solution of the homogeneous ODE.
By the linearity of $L$ we have
\begin{equation}
L(x_{p0}+x_{h0}) = L(x_{p0})+L(x_{h0})= r + 0 = r,
\end{equation}
and
\begin{equation}
[{\mathbf x}_{p0}+{\mathbf x}_{h0}](0) = {\mathbf x}_{p0}(0) + {\mathbf x}_{h0}(0) = {\mathbf 0} + {\mathbf x}_0 = {\mathbf x}_0.
\end{equation}
Thus $x_{p0}+x_{h0}$ is the solution to (\ref{IVP}).

The advantage of this particular decomposition is seen, when the {\em Laplace-transform} is brought into play and thus we are able to define the powerful tool of the {\em network function} of an LTI-circuit (cf.\ \cite{ChuaDesoerKuh}, Chap.\ 10, sect. 4.4).

\subsection{Transient and Steady-State}
A decidedly different decomposition is used, when we demand that the right-hand side $r$ of (\ref{IVP}) for some $\omega\in {\mathbf R}$ satisfies
\begin{equation}
r\in {\mathfrak C}_\omega,
\end{equation}
where ${\mathfrak C}_\omega\subset C^n(I)$ is defined by 
\begin{equation}
{\mathfrak C}_\omega := \{\, f:I\to {\mathbf R}, t\mapsto \alpha\cos(\omega t) + \beta\sin(\omega t)\,\vert\, \alpha,\beta\in {\mathbf R}\,\},
\end{equation}
and all the zeroes of the characteristic polynomial $p_L(\lambda)$ have negative $(< 0)$ real parts. Under these assumptions, there is a {\em bounded} function $x_{pr}$, which is a solution to the inhomogeneous ODE
\begin{equation}
\label{steadystate}
L(x_{pr}) = r
\end{equation}
and which is {\em uniquely determined already by} (\ref{steadystate}).
Furthermore we have
\begin{equation}
x_{pr}\in {\mathfrak C}_\omega.
\end{equation}
The function $x_{pr}$ is called the {\em steady-state solution} of (\ref{IVP}).
The corresponding homogeneous solution $x_{hr}$, which satisfies 
\begin{equation}
x_0 = x_{pr} + x_{hr}
\end{equation}
is referred to as the {\em transient solution}. 

The method to determine $x_{pr}$ is known as {\em sinusoidal steady-state} or {\em phasor analysis} (cf.\ \cite{ChuaDesoerKuh}, Chap.\ 9). While this method is considered by some as theoretically unsound (cf.\ e.g.\ \cite{Control}, Chap.\ 3.3., where it is called a \lq\lq rule of thumb\rq\rq), the work of Marten and Mathis \cite{MM1}, \cite{MM2}, \cite{MM3}, \cite{MM4}, \cite{MartenET}, as will be argued below, has shown that this is not the case.

\subsection{Steady-State response vs unbounded input}
While the above cited result on the existence and uniqueness of the steady-state solution can be generalized to inputs $r\in C_b(I)$, i.e.\ to arbitrary {\em bounded} continuous input functions \cite{MartenET}, it has become customary in Control Theory to consider unbounded input and speak of the resulting response, which usually is unbounded, too, as a \lq\lq steady-state\rq\rq, as well (cf.\ \cite{Control}, Chap.\ 3.1 and 3.4). We point out the fact, that this linguistic lapse can be corrected \cite{steady}, but we will not do so in this paper. Even worse, from this point onward, we will make the same use of the word \lq\lq steady-state\rq\rq\ and will call any solution $x_p$ of the equation $L(x_p) = r$ a {\em steady-state-solution}, keeping in mind that in general, while existence is guaranteed for any $r\in C(I)$ by theorem \ref{Walter}, uniqueness usually is not.

\subsection{Some more notation}
We close this section by fixing the notation for some more sets of input-functions. Let $\lambda, \omega \in {\f R}, m\in {\f N}\cup \{0\}$. With ${\f R}[t]$, we denote the $\f R$-linear space of polynomial functions. We define
\begin{align}
{\g E}_\lambda &:= \{\, t \mapsto \alpha e^{\lambda t}\,:\, \alpha \in {\f R}\,\},\\
{\g E}_{\lambda+j\omega} &:=\{\, t \mapsto \alpha \cos(\omega t) e^{\lambda t}+ \beta \sin(\omega t)e^{\lambda t}\,:\, \alpha,\beta \in {\f R}\,\},\\
{\g {EP}}_\lambda &:= \{\, t \mapsto p(t)\cdot e^{\lambda t}\,:\, p(t)\in {\f R}[t]\,\},\\
{\g {EP}}^m_\lambda &:= \{\, t \mapsto p(t)\cdot e^{\lambda t}\,:\, p(t)\in {\f R}[t], \deg(p(t)) \le m\,\},\\
{\g {CP}}_\omega &
:= 
\{\, t \mapsto p(t)\cdot \cos(\omega t) + q(t)\cdot \sin(\omega t)\,:\,\notag\\
&\phantom{:= \{\,xxxxxxxxxxxxxxxxxxx}\quad
          p(t),q(t)\in {\f R}[t]\,\},\\
{\g {CP}}^m_\omega &
:= 
\{\, t \mapsto p(t)\cdot \cos(\omega t) + q(t)\cdot \sin(\omega t)\,:\,\notag\\
&  p(t),q(t)\in {\f R}[t], \max\{\deg(p(t)), \deg(q(t))\}\le m\,\},\\
{\g {EP}}_{\lambda+j\omega} &:=\{\, t \mapsto p(t)\cos(\omega t) e^{\lambda t}+ q(t) \sin(\omega t)e^{\lambda t}\,:\,\notag\\
&\phantom{:= \{\,xxxxxxxxxxxxxxxxxxx}\quad
          p(t),q(t)\in {\f R}[t]\,\},\\
\intertext{and finally}
{\g {EP}}^m_{\lambda+j\omega} &:=\{\, t \mapsto p(t)\cos(\omega t) e^{\lambda t}+ q(t) \sin(\omega t)e^{\lambda t}\,:\,\notag\\
&  p(t),q(t)\in {\f R}[t], \max\{\deg(p(t)), \deg(q(t))\}\le m\,\}.
\end{align}

It is obvious that all of the above sets form $\f R$-linear spaces of real-valued functions -- ${\g {CP}}_\omega$ and ${\g {EP}}_{\lambda+j\omega}$ being infinite dimensional -- and we clearly have the inclusions
\begin{equation}
{\g C}_\omega \subset {\g {CP}}^m_\omega \subset {\g {CP}}_\omega, 
\end{equation}
\begin{equation}
{\g E}_{\lambda+j\omega}\subset {\g {EP}}^m_{\lambda+j\omega}\subset {\g {EP}}_{\lambda+j\omega} 
\end{equation}
and 
\begin{align}
{\g C}_\omega &= {\g E}_{j\omega},\\
{\g {CP}}^m_\omega &= {\g {EP}}^m_{j\omega},\\
{\g {CP}}_\omega &= {\g {EP}}_{j\omega}.
\end{align}
For $\lambda = \omega = 0$ we have, that ${\g {EP}}$ is equal to the space ${\f R}[t]$ of polynomial functions.

In addition, we have that each of these spaces is closed under differentiation, thus $L$ induces a linear operator on each of them. The following sections are dedicated to a closer study of the action of $L$ in each of the above cases.

\section{Complex Structures on the above function spaces}
Earlier we have stressed the fact, that all of the above sets of functions are real linear spaces. Marten and Mathis have shown in the eighties that ${\g C}_\omega$ carries the structure of a {\em complex} linear space, as well, and from this were able to build up the AC-calculus for functions in ${\g C}_\omega$. In this and the next sections, step-by-step we will generalize this result to the space ${\g {EP}}_{\lambda+j\omega}$, thus setting the foundation for the \lq\lq full\rq\rq\ AC-calculus. We will further look at how differentiation and, consequently, the action of the operator $L$ fit together with this complex structure.

\subsection{AC-Calculus for Functions in ${\g C}_\omega$ vs. Phasor Analysis}
In this section we closely follow the line of thought set up in \cite{MM1}, \cite{MM2}, \cite{MM3}, \cite{MM4} and \cite{MartenET}.

Let $x(t)= \alpha\cos(\omega t) + \beta\sin(\omega t)\in {\g C}_\omega$ and 
$a+jb\in {\f C}$. Then we can define a multiplication $\odot$ with complex scalars on ${\g C}_\omega$ by setting
\begin{eqnarray*}
j\odot x(t) & = & j \odot (\alpha\cos(\omega t) + \beta\sin(\omega t))\\
& := & - \alpha\sin(\omega t) + \beta\cos(\omega t)
\end{eqnarray*}
and continuing this rule by
\begin{eqnarray}
&& (a+jb) \odot x(t)\notag\\
& = & (a+jb) \odot (\alpha\cos(\omega t) + \beta\sin(\omega t))\notag\\
& := & (a\alpha+b\beta)\cdot \cos(\omega t) + (a\beta - b\alpha)\cdot\sin(\omega t).
\end{eqnarray}
With this multiplication at hand, ${\g C}_\omega$ becomes a $1$-di\-men\-sio\-nal {\em complex} linear space of functions. If we fix as a basis the function $t \mapsto \cos(\omega t)$, then we have
\begin{equation}
\alpha\cos(\omega t) + \beta\sin(\omega t) = (\alpha-j\beta)\odot \cos(\omega t).
\end{equation}
Let us point out the fact that, if we consider a {\em sinusoid} 
\begin{equation}
t\mapsto z(t) = A\cos(\omega t + \Phi)
\end{equation}
with $A, \omega,\Phi\in {\f R}$, then due to the addition rule for the cosine, we have
\begin{eqnarray}
A\cos(\omega t + \Phi) & = & A\cdot\left(\cos(\Phi)\cos(\omega t)-\sin(\Phi)\sin(\omega t)\right)\notag\\
& = & [A\cdot(\cos(\Phi)+j\sin(\Phi))]\odot\cos(\omega t)\notag\\
& = & (Ae^{j\Phi})\odot\cos(\omega t).
\end{eqnarray}
Thus {\em the phasor associated to the sinusoid} $z(t)$ {\em is the complex scalar with which the basis-function} $\cos(\omega t)$ {\em has to be multiplied according to} $\odot$ {\em to get} $z(t)$. 

\smallskip
We now look at the effects of differentiation on functions in ${\g C}_\omega$. For $x(t)$ as above we find
\begin{eqnarray}
\frac{d}{dt} x(t) & = & \frac{d}{dt} (\alpha\cos(\omega t) + \beta\sin(\omega(t)))\notag\\
& = & \omega\beta\cos(\omega t) - \omega\alpha\sin(\omega(t))\notag\\
& = & (j\omega)\odot(\alpha\cos(\omega t) + \beta\sin(\omega(t)))\\
& = & (j\omega)\odot x(t).
\end{eqnarray}
Thus {\em differentiating a function} $x(t)\in {\g C}_\omega$ {\em is the same as multiplying} $x(t)$ {\em with the complex scalar} $j\omega$ {\em according to} $\odot$.

\smallskip
Now as an easy consequence, we have for all functions $x\in {\g C}_\omega$
\begin{equation}
L(x) = p_L(j\omega)\odot x
\end{equation}
where $p_L$ again denotes the characteristic polynomial of the operator $L$. Thus we get

\smallskip
\theorem
Let $r\in {\g C}_\omega$. If $j\omega$ is not a zero of the characteristic polynomial $p_L$ of the differential operator $L$, then the function
\begin{equation*}
x = \frac{1}{p_L(j\omega)}\odot r
\end{equation*}
is a particular solution to the inhomogeneous linear ODE
\begin{equation*}
L(x) = r.
\end{equation*}
\endproof

Thus -- if there had been any doubts -- phasor analysis now is completely rehabilitated.

\subsection{A slight generalization}
We are now looking at the space ${\g E}_{\lambda + j\omega}$. A typical function $x$ in this set is given by
\begin{equation}
x(t) = \alpha \cos(\omega t) e^{\lambda t}+ \beta \sin(\omega t)e^{\lambda t}
\end{equation}
with $\alpha, \beta\in {\f R}$.

Here again, we can endow ${\g E}_{\lambda + j\omega}$ with a complex structure
by setting
\begin{eqnarray}
&&(a+jb)\odot x(t)\notag\\
& = & (a+jb)\odot(\alpha \cos(\omega t) e^{\lambda t}+ \beta \sin(\omega t)e^{\lambda t})\notag\\
& := & (a\alpha + b\beta)\cos(\omega t) e^{\lambda t}+ (a\beta -b\alpha)
sin(\omega t)e^{\lambda t}.
\end{eqnarray}
${\g E}_{\lambda + j\omega}$ thus becomes a $1$-dimensional complex space, as well. 

The differential operator $\frac{d}{dt}$ acts on ${\g E}_{\lambda + j\omega}$
via
\begin{eqnarray}
&& \frac{d}{dt} x(t)\notag\\ 
& = & \frac{d}{dt}(\alpha \cos(\omega t) e^{\lambda t}+ \beta \sin(\omega t)e^{\lambda t})\notag\\
& = &
(\alpha\lambda + \beta\omega)\cos(\omega t) e^{\lambda t} + (\beta\lambda - \alpha\omega) \sin(\omega t)e^{\lambda t} \\
& = & (\lambda + j\omega) \odot (\alpha \cos(\omega t) e^{\lambda t}+ \beta \sin(\omega t)e^{\lambda t})\notag\\
& = & (\lambda + j\omega) \odot x(t).
\end{eqnarray}

With the same reasoning as in the above case, we thus attain the following result:

\smallskip
\theorem
Let $r\in {\g E}_{\lambda +j\omega}$. If $\lambda+j\omega$ is not a zero of the characteristic polynomial $p_L$ of the differential operator $L$, then the function
\begin{equation*}
x = \frac{1}{p_L(\lambda+j\omega)}\odot r
\end{equation*}
is a particular solution to the inhomogeneous linear ODE
\begin{equation*}
L(x) = r.
\end{equation*}
\endproof

\section{The full AC-Calculus}
The full AC-Calculus is a generalization to both phasor analysis, as described above, and the classic \lq\lq Ansatz according to the right hand side\rq\rq. Thus it relies heavily on linear algebra.

\subsection{Polynomial input and input from ${\g {EP}}^m_\lambda$}
Let us fix $m\in {\f N}$. We know, that ${\g {EP}}^m_\lambda$ is an $m+1$-dimensional real linear space of functions.
For the duration, we take as a fixed basis in ${\g {EP}}^m_\lambda$ the functions
$\{\, u_k\,:\, k\in \{\,0,\dots,m\,\}\}$, where we set
\begin{equation}
u_k :t \mapsto t^k\cdot e^{\lambda t}.
\end{equation} 
As we know from theorem \ref{Walter} those functions are linearly independent and, indeed, form a basis.

${\g {EP}}^m_\lambda$ is closed under differentiation, thus $\frac{d}{dt}$ gives a linear operator on ${\g {EP}}^m_\lambda$ which is determined by its action on the basis functions
\begin{eqnarray}
\frac{d}{dt} u_k 
& = & \frac{d}{dt} t^k\cdot e^{\lambda t}\notag\\
& = & k\cdot t^{k-1}\cdot e^{\lambda t} + \lambda t^k\cdot e^{\lambda t}\notag\\
& = & k u_{k-1} + \lambda u_k
\end{eqnarray}
for $k\ge 1$ and
\begin{equation}
\frac{d}{dt} u_0 = \frac{d}{dt} e^{\lambda t} = \lambda e^{\lambda t} = \lambda u_0.
\end{equation}
If we identify the function
\begin{equation*}
x(t) = \sum_{k=0}^m \alpha_k t^k e^{\lambda t} = \sum_{k=0}^m \alpha_k u_k
\end{equation*}
with the vector
\begin{equation}
\underline{\f x} = (\alpha_0,\alpha_1,\dots,\alpha_m)^\top,
\end{equation} 
then the action of $\frac{d}{dt}$ on ${\g {EP}}^m_\lambda$ is given by matrix multiplication as
\begin{equation}
\frac{d}{dt} x \ \simeq\ 
\begin{pmatrix}
\lambda & 1 & & & & &\\
& \lambda & 2 & & & &\\
& & \lambda & 3 & & &\\
& & & \ddots & \ddots & &\\
& & & & \lambda & m-1 &\\
& & & & &  \lambda & m\\
& & & & & & \lambda\\
\end{pmatrix}
\underline{\f x}
\end{equation}
where empty entries denote $0$.
From this it follows, that the action of the operator $L$ is given by an $(m+1)\times (m+1)$-matrix as well. The next theorem can be shown by induction.

\smallskip
\theorem
Let $p_L$ again denote the characteristic polynomial of $L$. Then
$$
L(x) \simeq (a_{kl})_{0\le k,l\le m}\ \underline{\f x}
$$
where the matrix $(a_{kl})_{0\le k,l\le m}$ is given by
\begin{equation}
\label{akl}
a_{kl} = 
\begin{cases}
0 & \text{for $k > l$,}\\
\\
p_L(\lambda) & \text{for $k = l$,}\\
\\
\binom{l}{k}p_L^{(l-k)}(\lambda) & \text{for $k < l$.}\\
\end{cases}
\end{equation}
A slight exercise in matrix computation shows that this matrix $(a_{kl})$ is invertible, iff $\lambda$ is not a zero of $p_L$. Thus we have

\smallskip
\theorem
\label{EP}
If $r\in {\g {EP}}^m_\lambda$ is given by
\begin{equation}
r(t) = \sum_{k=0}^m r_k t^k  e^{\lambda t},
\end{equation}
and if
$\lambda$ is not a zero of the characteristic polynomial $p_L$, then a particular solution $x$ to the equation $L(x) = r$ is given by
\begin{equation}
x(t) = \sum_{k=0}^m x_k t^k e^{\lambda t},
\end{equation} 
where 
\begin{equation}
x_m = \frac{1}{p_L(\lambda)}\cdot r_m
\end{equation}
and successively
\begin{equation}
x_{m-l} = \frac{1}{p_L(\lambda)}\cdot \left( r_{m-l} - \sum_{k=m-l+1}^m \binom{k}{m-l} p^{(k-m+l)}(\lambda)x_k\right)
\end{equation}
for $1\le l\le m$.
\endproof
Setting $\lambda = 0$ in the above discussion settles the case of polynomial input.

\subsection{Input from ${\g {CP}}^m_{\omega}$}
First we show how to make ${\g {CP}}^m_{\omega}$ into a complex linear space. For
\begin{equation}
x(t) = p(t)\cos(\omega t) + q(t)\sin(\omega t)
\end{equation}
with real polynomials $p(t), q(t) \in {\f R}[t]$ and $a,b\in {\f R}$ we set
\begin{eqnarray}
&&(a+jb)\odot x(t)\notag\\
& = & (a+jb)\odot (p(t)\cos(\omega t) + q(t)\sin(\omega t))\notag\\
& := & (a\cdot p(t) + b\cdot q(t))\cos(\omega t)\notag\\
&&\quad\quad\quad\quad\quad + (a\cdot q(t) - b\cdot p(t))\sin(\omega t).
\end{eqnarray}
While ${\g {CP}}^m_\omega$ is a $2(m+1)$-dimensional real linear space, with
this new scalar multiplication $\odot$ it becomes an $m+1$-dimensional complex linear space. Fixing as a $\f C$-basis the set of functions $\{\, v_k\, :\, 0\le k\le m\,\}$, where
\begin{equation}
v_k : t\mapsto t^k \cos(\omega t)
\end{equation}
we see that differentiation acts as a $\f C$-linear operator on $x\in {\g {CP}}^m_\omega$. The corresponding matrix is given by the formula
\begin{equation}
\frac{d}{dt} x \ \simeq\ 
\begin{pmatrix}
j\omega & 1 & & & & &\\
& j\omega & 2 & & & &\\
& & j\omega & 3 & & &\\
& & & \ddots & \ddots & &\\
& & & & j\omega & m-1 &\\
& & & & &  j\omega & m\\
& & & & & & j\omega\\
\end{pmatrix}
\underline{\f x}.
\end{equation}
In analogy with the above discussion we get

\smallskip
\theorem
\label{CP}
If $r\in {\g {CP}}^m_\omega$ is given by
\begin{eqnarray*}
r(t) & = & \sum_{k=0}^m  t^k  (\alpha_k \cos(\omega t) + \beta_k \sin(\omega t))\\
& = & \sum_{k=0}^m \gamma_k \odot t^k \cos(\omega t)
\end{eqnarray*}
with $\gamma_k = \alpha_k - j\beta_k$ for $0\le k\le m$,
and if $j\omega$ is not a zero of the characteristic polynomial $p_L$,
then a particular solution $x$ to the equation $L(x) = r$ is given by
\begin{equation}
x(t) = \sum_{k=0}^m x_k \odot t^k \cos(\omega t),
\end{equation} 
where 
\begin{equation}
x_m = \frac{1}{p_L(j \omega)}\gamma_m
\end{equation}
and successively
\begin{eqnarray*}
&& x_{m-l}\notag\\
& = & \frac{1}{p_L(j\omega)}\cdot \left( \gamma_{m-l} - \sum_{k=m-l+1}^m \binom{k}{m-l} p^{(k-m+l)}(j\omega)x_k\right)
\end{eqnarray*}
for $1\le l\le m$.
\endproof

\section{Conclusion}
Obviously the full AC-Calculus, as demonstrated above, extends to the space ${\g {EP}}_{\lambda+j\omega}$, too, since we can define a multiplication with complex scalars for these functions, as well. 

Furthermore in this paper, we have not yet regarded the question, what will happen, when we have resonance, i.e. when in theorem \ref{EP} the scalar $\lambda$ or in theorem \ref{CP} the scalar $j\omega$ is a multiple root of the characteristic polynomial $p_L$. Keen readers will have guessed that the spaces ${\g {CP}}_\omega$, ${\g {EP}}_\lambda$ and ${\g {EP}}_{\lambda+j\omega}$ are perfectly suited for handling these cases. In fact, the formulas given in theorems \ref{EP} and \ref{CP} can easily be adapted. We will leave this task for a later time. 

For now, we hope that we were able to convince the gentle reader that it is perfectly alright and algorithmically advantageous to consider sets of real functions as complex linear spaces.

\bibliographystyle{IEEE}

\newpage
\begin{biography}
[{\includegraphics[width=1in,height=1.25in,clip,keepaspectratio]{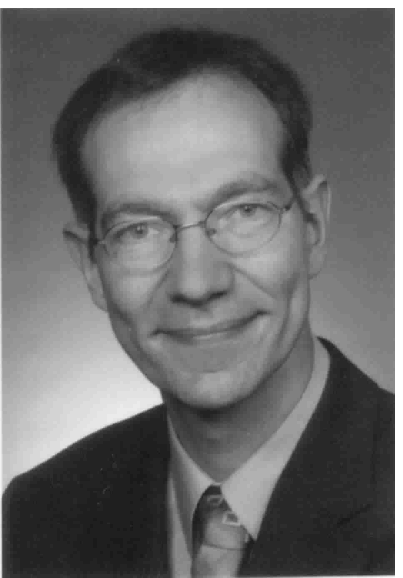}}]
{Eberhard H.-A.~Gerbracht}
received a Dipl.-Math.\ degree in mathematics, a Dipl.-Inform.\ degree in computer science, and a Ph.D. (Dr.\ rer.nat.) degree in mathematics from the Technical University Braunschweig, Germany, in 1990, 1993, and 1998, respectively.

From 1992 to 1997 he was a Research Fellow and Teaching Assistant at the Institute for Geometry at the TU Braunschweig. From 1997 to 2003 he was an Assistant Professor in the Department of Electrical Engineering and Information Technology at the TU Braunschweig. During that time he was also appointed lecturer for several courses on digital circuit design at the University of Applied Sciences Braunschweig/Wolfenb\"uttel, Germany. From 2001 to 2002 he was appointed lecturer for a two-semester course in linear circuit analysis at the TU Braunschweig. After a two-year stint as a mathematics and computer science teacher at a grammar school in Braunschweig and a vocational school in Gifhorn, Germany, he is currently working as an independent researcher in various areas of mathematics. His research interests include combinatorial and computer algebra, and their applications in geometry, calculus, and electrical engineering, C*-algebras, quantum computing, and the history of mathematics in the 19th and early 20th century.

Dr.~Gerbracht is a member of the German Mathematical Society (DMV), the German Society for Didactics of Mathematics (GDM), and the society ``Web Portal: History in Braunschweig - www.gibs.info''.
\end{biography}

\section*{Note added to the Electronic Version}
\small
In this electronic document, some small typographical errors of the printed version were corrected. This especially refers to the last formula given in Theorem \ref{CP}.

Furthermore, for the convenience of the reader an abstract, keywords, MSC classification, and a short CV according to IEEE standards have been added to the arXiv-version. 
\hfill (Sept.~17th,~2007)

\end{document}